\documentclass[12pt]{amsart}
\usepackage{verbatim}
\usepackage{amsmath}
\usepackage{amssymb}

\theoremstyle{plain}
\newtheorem*{thm}{Theorem}

                                                \newtheorem*{lem}{Lemma}

\makeatletter
 \def\@setcopyright{}
      \def\serieslogo@{}
      \makeatother

\begin{document}

\fboxrule 1
\fboxsep
\author[Boris Weisfeiler]{\framebox{Boris Weisfeiler}} 
   %\address{Department of Mathematics, 
%Pennsylvania State University,     
%University Park,  PA 16802}   
%\email{\\}

\title{On one class of unipotent subgroups of semisimple algebraic groups$\ ^1$}

%\date{\today}
\begin{abstract}This paper contains a complete proof of a fundamental theorem
on the normalizers of unipotent subgroups in semisimple algebraic groups. A
similar proof was given later by A. Borel and J. Tits (see also 
bibliographical remarks in the introduction to their paper).$^2$

\end{abstract}

\footnotetext[1]{Translated from Russian, Uspehi Mat. Nauk, {\bf 21}, 2(128)
(1966), 222-223}

\footnotetext[2]{A. Borel, J. Tits, \emph{El\'ements unipotents et
sous--groupes paraboliques de groupes r\'eductifs I}, Invent. Math. {\bf 12}
(1971), 95--104. }
 
%\thanks{} 
\maketitle

I.I. Pyatetskii--Shapiro suggested me to prove the following theorem which  in
his opinion would be of interest for the theory of discrete subgroups of
semisimple Lie groups.  
\begin{thm} Let $k$ be an arbitrary field,
$G$ a semisimple algebraic group defined over $k$, and $H$ a unipotent subgroup
of
$G$. If the nilpotent radical of the normalizer $N_G(H)$ of  $H$ in
$G$,  coincides with
$H$, then
$N_G(H)$  is a parabolic subgroup of $G$.
\end{thm}

\begin{lem} Let $N$ and $H$ be nilpotent groups, $H\subset N$. If
$N_N(H)=H$, then $N=H$.
\end{lem}

Let $N_0=N$, $N_{s+1}=[N_s,N]$ with
$N_q=1$. Then $N_q\subset H$. Further, if $N_{s+1}\subset H$, then $N_s\subset
H$ since $[N_s,H]\subset [N_s,N]= N_{s+1}\subset H$. Thus $N_s\subset H$ for any
$s$, which completes the proof.

\begin{proof}[Proof of the Theorem] First we prove that $N_G(H)$ contains a
maximal torus $T$ of the group $G$. For, let $N_G(H)=SH$ where $S$ is a
reductive group. Let us choose two Borel subgroups $B$ and $B'$ in $G$
satisfying the following conditions:
$B\cap N_G(H)$ and $B'\cap N_G(H)$ are Borel subgroups of $N_G(H)$, and
$B\cap B'\cap N_G(H)=\tilde TH$, where $\tilde T$ is a maximal torus of  $S$.
Let $N$ be the unipotent part of  $B\cap B'$. Then
$N_N(H)=N_G(H)\cap N=H$, and by lemma, $N=H$. But, it is well--known (\cite{T},
n$^o$ 2.16) that $B\cap B'=TN$ where $T$ is a maximal torus of $G$ and
$N$ is a normal subgroup of $TN$, whence $T\subset N_G(H)$.

Assume now that $T\subset S$, $T_1$ is the center of $S$,
$\Gamma=N_S(T)/T$, and $\Sigma$ and $\Sigma^+$ are the system of roots and the
system of positive roots of the group $G$ with respect to 
$T$, respectively. The group
$\Gamma$ is contained in the Weyl group $N_G(T)/T$ of  $G$, and therefore acts
on $\Sigma$. Since $T\subset N_G(H)$, $H$ is generated by the root subgroups
$P_\alpha,\,\alpha\in\Sigma^+$, which it contains (\cite {C}, exp 13, t. 1). If
$P_\alpha\subset H$ then $P_\alpha^\sigma\subset H$ for any$\sigma\in \Gamma$
(since $N_s(T)\subset S\subset N_G(H)$). Hence $\alpha^\sigma\in\Sigma^+$ for
any $\sigma\in \Gamma$. Let $N$ be a group generated by those
$P_\alpha,\,\,\alpha\in\Sigma^+$, for which
$\alpha^\sigma\in\Sigma^+$ for any $\sigma\in \Gamma$. Then $N\supset H$,
$N_N(H)=H$ and therefore $H=N$. It is known (\cite {S}, Prop. 3)
 that the group
$Z_G(T_1)N$ is parabolic and its nilpotent radical is $N$ (here $Z_G(T_1)$ is
the centralizer of $T_1$ in $G$). This completes the proof of the Theorem.
\end{proof}

The author takes the opportunity to express his gratitude to E.B. Vinberg and
I.I. Pyatetskii--Shapiro.

\bigskip

\noindent \emph{Remark added by E.B. Vinberg in translation.} The 
reference \cite {S} is 
not adequate, but the fact the author 
needs is quite trivial. One can proceed as  follows. The subgroup $P=Z_G(T_1)B$
is parabolic, with a maximal reductive subgroup  $Z_G(T_1)\supset S$ and some
unipotent radical  $U$. Since the semisimple parts of  $Z_G(T_1)$ and  $S$ have
the same rank, for any root $\alpha$ of  $Z_G(T_1)$ there is $\sigma\in \Gamma$
such that
$\alpha^\sigma\in \Sigma^-$. It follows  that $H\subset U$ and $N_U(H)=SH\cap
U=H$, whence $U=H$ and  $Z_G(T_1)=S$.


\begin{thebibliography}{99}
\bibitem {C} C. Chevalley, \emph {Classification des groupes de Lie
algebraiques}, t. 1--2,  E.N.S., Paris, 1956--1958

\bibitem {S} I. Satake, \emph{On the theory of reductive algebraic groups over
a perfect field}, J. Math. Soc. Japan, {\bf 15} (1963), 210--235

\bibitem {T} J. Tits, \emph{Groupes semi--simple isotropes}, in \emph{Colloq.
Th\'eorie des groupes Alg\'ebraiques} (Bruxelles, 1962), 137--147
\end{thebibliography}
\end{document}